\documentclass[draft]{amsart}

\usepackage{amssymb}
\usepackage{url}
\usepackage{amscd}
\usepackage{mathrsfs}
\usepackage{amsmath}
\usepackage{amsfonts}
\usepackage{amsfonts,enumerate}
\usepackage{pifont}
\usepackage{enumerate}
\usepackage{graphics}
\usepackage{verbatim}
\usepackage{amsthm}

\theoremstyle{plain}
\newtheorem{theorem}{Theorem}

\newtheorem{corollary}[theorem]{Corollary}
\newtheorem{lemma}[theorem]{Lemma}
\newtheorem{proposition}[theorem]{Proposition}
\newtheorem{definition}[theorem]{Definition}

\theoremstyle{remark}
\newtheorem{remark}{Remark}

\begin{document}


\title{Hardy spaces of operator-valued analytic functions}

\author{Zeqian Chen}

\address{Wuhan Institute of Physics and Mathematics,
Chinese Academy of Sciences\\ 30 West District, Xiao-Hong-Shan,
Wuchang, Wuhan 430071, China}

\email{zqchen@wipm.ac.cn}

\thanks{Partially supported by NSFC grant No.10775175}


\thanks{Hardy and BMOA spaces, noncommutative $L^p$
spaces, Lusin area integral, Littlewood-Paley $g$-functions}

\begin{abstract}
We are concerned with Hardy and BMO spaces of operator-valued
functions analytic in the unit disk of $\mathbb{C}.$ In the case of
the Hardy space, we involve the atomic decomposition since the usual
argument in the scalar setting is not suitable. Several properties
(the Garsia-norm equivalent theorem, Carleson measure, and so on) of
BMOA spaces are extended to the operator-valued setting. Then, the
operator-valued $\mathrm{H}^1$-BMOA duality theorem is proved.
Finally, by the $\mathrm{H}^1$-BMOA duality we present the Lusin
area integral and Littlewood-Paley $g$-function characterizations of
the operator-valued analytic Hardy space.
\end{abstract}

\maketitle


\newcommand\sfrac[2]{{#1/#2}}

\newcommand\cont{\operatorname{cont}}
\newcommand\diff{\operatorname{diff}}


\section{Introduction}
The classical Hardy and BMO spaces and $\mathrm{H}^1$-BMO duality
theorem play a crucial role in harmonic analysis (see for example,
\cite{G, K, S93}). The vector-valued analogue was studied by
Bourgain \cite{Bour} and Blasco \cite{B} in the case of the unit
disc. Recently, operator-valued (= quantum) harmonic analysis has
developed considerably (e.g., see \cite{JMX}). This is inspired by
the works on matrix-valued harmonic analysis (e.g., see \cite{NPTV,
PS} and references therein) and the recent development on the theory
of non-commutative martingales (see the survey paper by Xu
\cite{X}). The theory of operator-valued Hardy and BMO spaces in
$\mathbb{R}$ was well built by Mei \cite{M}. The goal of this paper
is to study the disk analogue of Mei's results, with an emphasis on
the `analytical' aspect.

The remainder of this paper is divided into five sections. In
Section 2 we present some preliminaries, including the $L^p$-spaces
of operator-valued measurable functions, the noncommutative
H\"{o}lder inequality, operator-valued analytic functions and some
properties. The operator-valued analytic Hardy space is defined by
the atomic decomposition in Section 3, since the usual argument in
the scalar setting is not suitable, as was pointed out in \cite{M}.
In section 4 we show the Garsia-norm equivalent theorem and Carleson
measure characterization of the operator-valued BMOA space. One of
main results in the paper is then proved in section 5: the
operator-valued $\mathrm{H}^1$-BMOA duality theorem. Section 6 is
devoted to the proof of another main result, the Lusin area integral
and Littlewood-Paley $g$-function characterizations of the
operator-valued analytic Hardy space. Our argument uses the atomic
decomposition and is distinct from Mei's method originating from
\cite{PX1}, in which the Lusin area integral is used for defining
the Hardy space. The techniques involved here is slightly simpler
and also suitable for use in obtaining the corresponding results
found there. We would like to point out that the atomic
decomposition of the predual of operator-valued BMO spaces has been
studied by Ricard \cite{R}.

In what follows, $C$ always denotes a constant, which may be
different in different places. For two nonnegative (possibly
infinite) quantities $X$ and $Y$ by $X \asymp Y$ we means that there
exists a constant $C$ such that $C^{-1} X \leq Y \leq C X.$ Any
notation and terminology not otherwise explained, are as used in
\cite{G} for complex harmonic analysis, and in \cite{T} for theory
of von Neumann algebras.

\section{Preliminaries}
\subsection{Operator-valued measurable functions}
Throughout this paper, $\mathcal{M}$ will always denote a von
Neumann algebra, and $\mathcal{M}_+$ its positive part. Recall that
a {\it trace} on $\mathcal{M}$ is a map $\tau: \mathcal{M}_+
\rightarrow [0, \infty ]$ satisfying:
\begin{enumerate}[{\rm (1)}]

\item $\tau (x + y) = \tau (x) + \tau (y)$ for
arbitrary $x , y \in \mathcal{M}_+;$

\item $\tau ( \lambda x ) = \lambda \tau (x)$ for any $\lambda
\in [0, \infty )$ and $x \in \mathcal{M}_+;$ and

\item $\tau (u^* u) = \tau (u u^*)$ for all $u \in \mathcal{M}.$

\end{enumerate}
$\tau$ is said to be {\it normal} if $\sup_{\gamma} \tau
(x_{\gamma}) = \tau ( \sup_{\gamma} x_{\gamma})$ for each bounded
increasing net $(x_{\gamma} )$ in $\mathcal{M}_+,$ {\it semifinite}
if for each $x \in \mathcal{M}_+ \backslash \{0\}$ there is a
non-zero $y \in \mathcal{M}_+$ such that $y \leq x$ and $\tau (y) <
\infty,$ and {\it faithful} if for each $x \in \mathcal{M}_+
\backslash \{0\},$ $\tau (x) > 0.$ A von Neumann algebra
$\mathcal{M}$ is called {\it semifinite} if it admits a normal
semifinite faithful trace $\tau,$ which we assume in the remainder
of this paper.

Denote by $\mathcal{S}_+$ the set of all $x \in \mathcal{M}_+$ such
that $\tau (\mathrm{supp} x) < \infty,$ where $\mathrm{supp} x$ is
the support of $x$ which is defined as the least projection $p$ in
$\mathcal{M}$ so that $p x = x$ or equivalently $x p = x.$ Let
$\mathcal{S}$ be the linear span of $\mathcal{S}_+.$ Then
$\mathcal{S}$ is a $\ast$-subalgebra of $\mathcal{M}$ which is
$w^*$-dense in $\mathcal{M}.$ Moreover, for each $0 < p < \infty,$
$x \in \mathcal{S}$ implies $| x |^p \in \mathcal{S}_+$ (and so
$\tau (| x |^p)< \infty$), where $|x| = (x^* x )^{1/2}$ is the
modulus of $x.$ Now we define $\| x \|_p = \left [\tau (| x |^p)
\right ]^{1/p} $ for all $x \in \mathcal{S}.$ One can show that
$\|\cdot \|_p$ is a norm on $\mathcal{S}$ if $1 \leq p < \infty,$
and a quasi-norm (more precisely, $p$-norm) if $0< p < 1.$ The
completion of $( \mathcal{S}, \|\cdot \|_p )$ is denoted by $L^p
(\mathcal{M}, \tau ).$ This is the non-commutative $L^p$-space
associated with $(\mathcal{M}, \tau ).$ For convenience, we set
$L^{\infty} (\mathcal{M}, \tau ) = \mathcal{M}$ equipped with the
operator norm. The trace $\tau$ can be extended to a linear
functional on $\mathcal{S},$ still denoted by $\tau.$ Since $| \tau
(x) | \leq \| x \|_1$ for all $x \in \mathcal{S},$ $\tau$ extends to
a continuous functional on $L^1 (\mathcal{M}, \tau ).$

Let $0 < \gamma, p, q \leq \infty$ be such that $1/ \gamma = 1/p +
1/q.$ If $x \in L^p (\mathcal{M}, \tau ), y \in L^q (\mathcal{M},
\tau )$ then $ xy \in L^{\gamma} (\mathcal{M}, \tau )$
and\begin{equation}\label{2.1}\| x y \|_{\gamma} \leq \| x \|_p \| y
\|_q.\end{equation}In particular, if $\gamma = 1,$ $| \tau (x y) |
\leq \| x y \|_1 \leq \| x \|_p \| y \|_q$ for arbitrary $x \in L^p
(\mathcal{M}, \tau )$ and $y \in L^q (\mathcal{M}, \tau ).$ This
defines a natural duality between $L^p (\mathcal{M}, \tau )$ and
$L^q (\mathcal{M}, \tau ): \langle x, y \rangle = \tau (x y).$ For
any $1 \leq p < \infty$ we have $L^p (\mathcal{M}, \tau )^* = L^q
(\mathcal{M}, \tau )$ isometrically. Thus, $L^1 (\mathcal{M}, \tau
)$ is the predual $\mathcal{M}_*$ of $\mathcal{M},$ and $L^p
(\mathcal{M}, \tau )$ is reflexive for $1 < p < \infty.$ (For the
theory of non-commutative $L^p$-spaces, see the survey paper by
Pisier and Xu \cite{PX2} and references therein).

Let $(\Omega,\mathcal{F},\mu)$ be a $\sigma$-finite measure space. A
$\mathcal{S}$-valued function $\varphi$ on $\Omega$ is said to be
{\it simple} if it can be written as
\begin{equation}\label{2.2}\varphi = \sum_j a_j \chi_{F_j},
\end{equation}where $a_j \in \mathcal{S}$ and $\{ F_j \}$
is a finite set of measurable disjoint subsets of $\Omega$ with $\mu
(F_j) < \infty.$ Denote by $S(\Omega, \mathcal{S})$ the set of all
simple $\mathcal{S}$-valued functions on $\Omega.$ For such a simple
function $\varphi$ we define\begin{equation}\label{2.3}\| \varphi
\|_{L^p_c } = \Big \| \Big ( \sum_j a^*_j a_j \mu (F_j ) \Big
)^{\frac{1}{2}} \Big \|_p = \left \| \left ( \int_{\Omega} | \varphi
|^2 d \mu \right )^{\frac{1}{2}} \right \|_p \end{equation}and $\|
\varphi \|_{L^p_r } = \| \varphi^* \|_{L^p_c },$ respectively. As
shown in \cite{M}, each $\varphi$ can be regarded as an element $T
(\varphi )$ in $\mathcal{M} \otimes \mathcal{B} (L^2 (\Omega))$
(where $\mathcal{B}( L^2 (\Omega))$ is the space of all bounded
operators on $L^2 (\Omega)$ with the usual trace $\mathrm{Tr}$)
and\begin{equation}\label{2.4}\| \varphi \|_{L^p_c } = \| T (\varphi
) \|_{L^p (\mathcal{M} \otimes \mathcal{B} (L^2
(\Omega)))},\end{equation}(and, $\| \varphi \|_{L^p_r } = \|
\varphi^* \|_{L^p_c }$). This concludes that $\| \cdot \|_{L^p_c }$
(and, $\| \cdot \|_{L^p_r }$) are norms on $S(\Omega, \mathcal{S})$
if $1 \leq p \leq \infty$ and $p$-norms if $0 < p <1.$ For $0 < p <
\infty$ (or, $p = \infty$) the completions of $S(\Omega,
\mathcal{S})$ in $\| \cdot \|_{L^p_c }$ and $\| \cdot \|_{L^p_r }$
(or, in weak$^*$-operator topology) are denoted by $L^p
(\mathcal{M}, L^2_c (\Omega))$ and $L^p (\mathcal{M}, L^2_r
(\Omega)),$ respectively.

\begin{lemma}\label{le2.1} (Proposition 1.1 in \cite{M}) Let $0 <
\gamma, p, q \leq \infty$ be such that $1/ \gamma = 1/p + 1/q.$ If
$f \in L^p (\mathcal{M}, L^2_c (\Omega)), g \in L^q (\mathcal{M},
L^2_c (\Omega)),$ then\begin{equation}\label{2.5} \langle f, g
\rangle = \int_{\Omega} f^* g d \mu \end{equation}exists as an
element in $L^{\gamma} (\mathcal{M}, \tau )$
and\begin{equation}\label{2.6}\| \langle f, g \rangle \|_{\gamma}
\leq \| f \|_{L^p_c } \| g \|_{L^q_c }. \end{equation}A similar
statement also holds for $L^p (\mathcal{M}, L^2_r (\Omega)).$
\end{lemma}

For $1 \leq p < \infty$ with $q = p / (p-1),$ by Lemma \ref{le2.1}
we have\begin{equation*}\label{}L^p (\mathcal{M}, L^2_c (\Omega))^*
= L^q (\mathcal{M}, L^2_c (\Omega))
\end{equation*}isometrically via the
anti-duality\begin{equation}\label{2.7}(f, g) = \tau ( \langle f, g
\rangle ) = \tau \left ( \int_{\Omega} f^* g d \mu \right )
\end{equation}for $f \in L^p (\mathcal{M}, L^2_c (\Omega))$
and $g \in L^q (\mathcal{M}, L^2_c (\Omega)).$
Similarly,\begin{equation*}L^p (\mathcal{M}, L^2_r (\Omega))^* = L^q
(\mathcal{M}, L^2_r (\Omega)).
\end{equation*}(For details see \cite{M}.)

Also, by the convexity of the operator-valued function $x \to
|x|^2,$ we have\begin{equation}\label{2.8}\Big |\int_{\Omega} f g d
\mu \Big |^2 \leq \int_{\Omega} |f|^2 d \mu \int_{\Omega} |g|^2 d
\mu \end{equation}for every $f \in L^p (\mathcal{M}, L^2_c
(\Omega))$ and $g \in L^2 (\Omega, \mu).$ (e.g., see (1.13) in
\cite{M}.)

\

\begin{lemma}\label{le2.2} If $f \in L^1 (\mathcal{M}, L^2_c
(\Omega)),$ then \begin{equation}\label{2.9}\left ( \int_{\Omega} |
\tau (f)|^2 d \mu \right )^{1/2} \leq \| f \|_{L^1_c}.
\end{equation}Consequently, $L^1 (\mathcal{M}, L^2_c
(\Omega)) \subset L^2 ( \Omega, L^1 (\mathcal{M})).$
\end{lemma}

\begin{proof}
Since$$\| g \|_{L^{\infty}_c} = \left \| \int_{\Omega } |g |^2 d \mu
\right \|^{1/2}_{\infty} \leq \left ( \int_{\Omega} \| g
\|^2_{\infty} d \mu \right )^{1/2} = \| g \|_{L^2 (\Omega ,
L^{\infty} (\mathcal{M}))},$$we conclude from \eqref{2.7} that
\begin{equation*}\label{}\begin{split}\left ( \int_{\Omega} |
\tau (f)|^2 d \mu \right )^{1/2} & \leq \| f \|_{L^2
(\Omega , L^1 (\mathcal{M}))}\\
& = \sup_{g \in L^2 (\Omega , L^{\infty} (\mathcal{M})), \| g \|_2
\leq 1} \left | \tau \left ( \int_{\Omega} f g^* d \mu \right )
\right | \\
& \leq \sup_{g \in L^{\infty} (\mathcal{M},L^2_c (\Omega )), \| g
\|_{L^{\infty}_c} \leq 1} \left | \tau \left ( \int_{\Omega} f g^* d
\mu \right ) \right | \\
& = \| f \|_{L^1_c}.\end{split}\end{equation*}This completes the
proof.
\end{proof}

For $0 < p \leq 2$ we set
\begin{equation}\label{2.10}L^p (\mathcal{M}, L^2_{cr}
(\Omega)) = L^p (\mathcal{M}, L^2_c (\Omega)) + L^p (\mathcal{M},
L^2_r (\Omega)),\end{equation}equipped with the sum
norm\begin{equation*}\| f \|_{L^p_{cr}} = \inf \left \{ \| g
\|_{L^p_c} + \| h \|_{L^p_r} \right \},\end{equation*}where the
infimum is taken over all $g \in L^p (\mathcal{M}, L^2_c (\Omega))$
and $h \in L^p (\mathcal{M}, L^2_r (\Omega))$ such that $f = g + h.$
For $2 < p \leq \infty,$\begin{equation}\label{2.11}L^p
(\mathcal{M}, L^2_{cr} (\Omega)) = L^p (\mathcal{M}, L^2_c (\Omega))
\cap L^p (\mathcal{M}, L^2_r (\Omega)),
\end{equation}equipped with the maximum norm$$\| f
\|_{L^p_{cr}} = \max \left \{ \| f \|_{L^p_c},  \| f \|_{L^p_r}
\right \}.$$Then, for $1 \leq p < \infty$ with $q = p / (p-1)$ we
have\begin{equation*}L^p (\mathcal{M}, L^2_{cr} (\Omega))^* = L^q
(\mathcal{M}, L^2_{cr} (\Omega))
\end{equation*}isometrically via the anti-duality \eqref{2.7}.

\subsection{Operator-valued analytic functions}
Let $\mathbb{D} = \{z \in \mathbb{C}: |z| < 1 \}$ be the unit disc
in the complex plane $\mathbb{C}$ and let $\mathbb{T} = \{z \in
\mathbb{C}: |z| = 1 \}$ be the unit circle. $d m = d \theta / 2 \pi$
will denote the normalized Lebesgue measure on $\mathbb{T}.$ The
kernel\begin{equation}\label{2.12}P_z(w) = \frac{1- |z|^2}{|1 -
\bar{z} w |^2}~~~(z \in \mathbb{D}, w \in \mathbb{D} \cup
\mathbb{T}) \end{equation}is called the Poisson kernel in
$\mathbb{D}.$ By \eqref{2.6} we can define the Poisson integral
$P[f]$ of a function $f \in L^p (\mathcal{M}, L^2_c ( \mathbb{T}) )$
or $L^p (\mathcal{M}, L^2_r ( \mathbb{T}) )$
by\begin{equation}\label{2.13}P[f] (z) = \int_{\mathbb{T}} P_z(t)
f(t) d m (t) \in L^p (\mathcal{M}) \end{equation}for $z \in
\mathbb{D}.$ For simplicity, we still denote $P[f]$ by $f.$
Similarly, we define the Cauchy integral $\mathfrak{C} (f)$ of a
function $f \in L^p (\mathcal{M}, L^2_c ( \mathbb{T}) )$ or $L^p
(\mathcal{M}, L^2_r ( \mathbb{T}) )$ by\begin{equation}\label{2.14}
\mathfrak{C} (f) (z) = \int_{\mathbb{T}} \frac{ f (t)}{1 - \bar{t}
z} d m (t) \in L^p (\mathcal{M})
\end{equation}for $z \in \mathbb{D}.$

We let $\mathrm{Aut} ( \mathbb{D})$ be the M\"{o}bius group of all
automorphisms of $\mathbb{D}.$ Every M\"{o}bius transformation can
be written as\begin{equation*}\psi (z) = e^{i \theta}\frac{z -
z_0}{1 - \bar{z}_0 z },~~(z \in \mathbb{D})
\end{equation*}with $\theta$ real and $|\bar{z}_0| < 1.$

Let $1 \leq p \leq \infty.$ Recall that $f: \mathbb{D} \rightarrow
L^p ( \mathcal{M} )$ is said to be analytic if $f$ is the sum of a
power series in the $L^p ( \mathcal{M} )$-norm, that is, $ f ( z ) =
\sum^{\infty}_{n=0} a_n z^n$ for every $z \in \mathbb{D},$ where
$a_n \in L^p ( \mathcal{M} ).$ The class of all such functions is
denoted by $\mathcal{H} (\mathbb{D}, L^p ( \mathcal{M} ) ).$ If $f
\in \mathcal{H} (\mathbb{D}, L^p ( \mathcal{M} ) )$ then all order
partial derivatives of $f$ exist and belong still to $\mathcal{H}
(\mathbb{D}, L^p ( \mathcal{M} ) ).$ By Lemma \ref{le2.1}, one has
$\mathfrak{C} (f) \in \mathcal{H} (\mathbb{D}, L^p ( \mathcal{M} )
)$ for any $f$ in $L^p (\mathcal{M}, L^2_c ( \mathbb{T}) )$ or $L^p
(\mathcal{M}, L^2_r ( \mathbb{T}) ).$ We
set\begin{equation*}\mathcal{H}^p (\mathcal{M}, L^2_c ( \mathbb{T})
) =\{ f \in L^p (\mathcal{M}, L^2_c ( \mathbb{S}) ): P[f] \in
\mathcal{H} (\mathbb{D}, L^p ( \mathcal{M} ) ) \}.
\end{equation*}Similarly, we have $\mathcal{H}^p (\mathcal{M}, L^2_r ( \mathbb{T})
)$ and $\mathcal{H}^p (\mathcal{M}, L^2_{cr} ( \mathbb{T}) ).$

For any $f$ in $L^p (\mathcal{M}, L^2_c ( \mathbb{T}) ),$ the
gradient $\nabla f (z)$ is the $L^p ( \mathcal{M} )$-valued vector
$(\partial f / \partial x, \partial f / \partial y)$ and
\begin{equation*}\left | \nabla f (z) \right |^2 =
\left | \frac{\partial f}{\partial x} \right |^2 + \left |
\frac{\partial f}{\partial y} \right |^2 \in L^{p/2} ( \mathcal{M}
)\end{equation*}with $z = x + i y.$ In this notation we
have\begin{equation}\label{2.15} \left | \nabla f (z) \right |^2 = 2
| f' (z) |^2 \end{equation}if $f$ is analytic.

\begin{lemma}\label{le2.3}
If $f \in L^1 (\mathcal{M}, L^2_c (\mathbb{T})),$
then\begin{equation}\label{2.16} \int_{\mathbb{T}} |f - f(0) |^2 d m
= \frac{1}{\pi} \int_{\mathbb{D}} | \nabla f (z) |^2 \log
\frac{1}{|z|} d x dy. \end{equation} Consequently, if $f \in L^1
(\mathcal{M}, L^2_c (\mathbb{T}))$ then\begin{equation}\label{2.17}
\int_{\mathbb{T}} |f - f(w) |^2 P_w d m \asymp \int_{\mathbb{D}} P_w
(z) | \nabla f (z) |^2 (1-|z|^2) d x dy, \end{equation}for every $w
\in \mathbb{D}.$
\end{lemma}

\begin{proof}
By Lemma 3.1 in \cite{G} we
have\begin{equation*}\int_{\mathbb{T}} \overline{(\phi - \phi(0))}
(\varphi - \varphi(0) )  d m = \frac{1}{\pi} \int_{\mathbb{D}}
\overline{\nabla \phi (z)} \nabla \varphi (z) \log \frac{1}{|z|} d x
dy\end{equation*}for all $\phi, \varphi \in L^2 (\mathbb{T}).$ Then,
for any $\mathcal{S}_{\mathcal{M}}$-valued simple function $f =
\sum_j a_j \chi_{F_j}$ on $\mathbb{T},$ one has
that\begin{equation*}\begin{split}\int_{\mathbb{T}} |f - f(0) |^2 d
m & = \int_{\mathbb{T}} \Big | \sum_j a_j (\chi_{F_j} -
\chi_{F_j}(0)) \Big |^2 d m \\
& = \sum_{j,k} a^*_j a_k \int_{\mathbb{T}} \overline{(\chi_{F_j} -
\chi_{F_j}(0))} (\chi_{F_k} - \chi_{F_k}(0) )  d m\\
& = \frac{1}{\pi}  \sum_{j,k} a^*_j a_k \int_{\mathbb{D}}
\overline{\nabla \chi_{F_j} (z)} \nabla \chi_{F_k} (z) \log
\frac{1}{|z|} d x dy\\
& = \frac{1}{\pi} \int_{\mathbb{D}} | \nabla f (z) |^2 \log
\frac{1}{|z|} d x dy.\end{split}\end{equation*}Since the set of
$\mathcal{S}_{\mathcal{M}}$-valued simple functions is dense in $L^1
(\mathcal{M}, L^2_c (\mathbb{T}))$ and $\tau$ is faithful, it is
concluded that \eqref{2.16} holds for all $ f \in L^1 (\mathcal{M},
L^2_c (\mathbb{T})).$

To prove \eqref{2.17}, we first prove the case of $w=0,$ that
is,\begin{equation}\label{2.18} \int_{\mathbb{T}} |f - f(0) |^2 d m
\asymp \int_{\mathbb{D}} | \nabla f (z) |^2 (1-|z|^2) d x dy.
\end{equation}To this end, we note that $1 - |z|^2 \leq 2
\log (1/|z|), |z| < 1.$ Hence, by \eqref{2.16} one has
that\begin{equation*} \int_{\mathbb{D}} | \nabla f (z) |^2 (1-|z|^2)
d x dy \leq 2 \pi \int_{\mathbb{T}} |f - f(0) |^2 d m.
\end{equation*}To prove the reverse inequality, suppose that the
integral on the left hand is finite and denoted by $A.$ For $|z| >
1/4,$ we have the reverse inequality $\log (1/ |z|) \leq C
(1-|z|^2).$ This yields that\begin{equation*} \int_{1/4 < |z| < 1} |
\nabla f (z) |^2 \log \frac{1}{|z|} d x dy \leq C A.
\end{equation*}Also, let $\mathbb{H}$ be the Hilbert space on which
$\mathcal{M}$ acts. For $|z| < 1/4,$ the subharmonicity of $z
\rightarrow \| \nabla f(h) (z) \|^2$ for each $h \in \mathbb{H}$
gives that\begin{equation*}\begin{split}\langle | \nabla f (z) |^2
(h), h
\rangle & = \| \nabla f(h) (z) \|^2 \\
& \leq \frac{16}{\pi} \int_{|w - z| < 1/4} \| \nabla f(h) (w) \|^2 d
x d y\\
& \leq \frac{32}{\pi} \int_{|w| < 1/2} \| \nabla f(h) (w) \|^2 (1-
|w|^2) d x d y\\
& \leq \frac{32}{\pi} \langle A(h), h
\rangle.\end{split}\end{equation*}Hence,\begin{equation*} \int_{|z|
\leq 1/4 } | \nabla f (z) |^2 \log \frac{1}{|z|} d x dy \leq
\frac{32}{\pi} A \int_{|z| \leq 1/4 } \log \frac{1}{|z|} d x dy \leq
C¡¡A.\end{equation*}Using \eqref{2.16} we conclude the proof of
\eqref{2.18}.

Now, fix $w \in \mathbb{D}.$ The
identity\begin{equation}\label{2.19} \int_{\mathbb{T}} |f - f(w) |^2
P_w d m = \frac{1}{\pi} \int_{\mathbb{D}} | \nabla f (z) |^2 \log
\left | \frac{1- \bar{w} z}{z- \bar{w}} \right | d x dy
\end{equation}has the same proof as \eqref{2.16}. Using
the identity\begin{equation*}1- \left | \frac{z-w}{1- \bar{w}z}
\right |^2 = (1-|z|^2)P_w (z),\end{equation*}we similarly obtain
\eqref{2.17}. This completes the proof.
\end{proof}

For every $\alpha > 1$ and $t \in \mathbb{T},$ let
\begin{equation*}\Gamma_{\alpha}( t ) = \{ z \in
\mathbb{D}: | t - z| < \alpha (1- |z|) \}.
\end{equation*}In what follows, we fix $\alpha > 1.$

\begin{lemma}\label{le2.4}
If $f \in L^1 (\mathcal{M}, L^2_c (\mathbb{T})),$
then\begin{equation}\label{2.20} \int_{\mathbb{T}} \Big (
\int_{\Gamma_{\alpha} (t)} | \nabla f (z) |^2 d x d y \Big ) d m
\asymp \int_{\mathbb{T}} | f - f(0) |^2 d m .
\end{equation}
\end{lemma}

\begin{proof}
Let $\mathbb{H}$ be the Hilbert space on which $\mathcal{M}$ acts.
Note that as $|z| \to 1,$ $m \big \{ t \in \mathbb{T}: z \in
\Gamma_{\alpha}( t ) \big \} \asymp 1- |z|^2.$ By Lemma \ref{le2.3}
one has that\begin{equation*}\label{}\begin{split} \Big \langle
\int_{\mathbb{T}} \Big ( \int_{\Gamma_{\alpha} (t)} &
 | \nabla f (z) |^2  d x d y \Big ) d
m (h), h \Big \rangle \\
& = \int_{\mathbb{D}} \big \langle | \nabla f (z) |^2 (h), h \big
\rangle m \left \{ t \in \mathbb{D}: z \in
\Gamma_{\alpha}( t ) \right \} d x d y\\
& \asymp  \int_{\mathbb{D}} \big \langle | \nabla f (z) |^2
(h), h \big \rangle (1-|z|^2) d x d y \\
& \asymp \Big \langle \int_{\mathbb{T}} | f - f(0)|^2  d m (h), h
\Big \rangle\end{split}\end{equation*}for all $h \in \mathbb{H}.$
\end{proof}

\section{Operator-valued Hardy space}

Throughout, we always denote by $I$ a subarc of $\mathbb{T}$ and
$|I| = m(I).$

\begin{definition}\label{df3.1}
$a \in L^1 (\mathcal{M}, L^2_c (
\mathbb{T}) )$ is said to be an $\mathcal{M}_c$-atom, if
\begin{enumerate}[{\rm (i)}]

\item  $a$ is supported in a subarc $I$ of $\mathbb{T},$

\item  $\int_I a d \sigma =0,$ and

\item  $ \| a \|_{L^1_c} \leq |I|^{- 1/2}.$
\end{enumerate}
\end{definition}

By Lemma \ref{le2.2}, one concludes that an $\mathcal{M}_c$-atom is
also an $L^1 (\mathcal{M})$-valued $2$-atom on $\mathbb{T}.$ Hence,
for each $\mathcal{M}_c$-atom $a$ one has\begin{equation}\label{3.1}
\| a \|_{L^1(\mathbb{T}, L^1(\mathcal{M}))} \leq |I |^{-1/2} \| a
\|_{L^2(\mathbb{T}, L^1(\mathcal{M}))} \leq 1. \end{equation}Then,
we define $\mathrm{H}^1_c (\mathbb{T}, \mathcal{M} )$ as the space
of all $f \in L^1 ( \mathbb{T}, L^1 (\mathcal{M} ))$ which admit a
decomposition\begin{equation*}f = \sum_k \lambda_k a_k ,
\end{equation*}where for each $k,$ $a_k$ is an $\mathcal{M}_c$-atom
or an element in $L^1 (\mathcal{M})$ of norm less than $1,$ and
$\lambda_k \in \mathbb{C}$ so that $\sum_k |\lambda_k| < \infty.$ We
equip this space with the norm\begin{equation*}\| f
\|_{\mathrm{H}^1_c} = \inf \Big \{ \sum_k | \lambda_k |: f = \sum_k
\lambda_k a_k \Big \},\end{equation*}where the infimum is taken over
all decompositions of $f$ described above.

\begin{proposition}\label{prop3.1}
Let $\mathrm{H}^1 (\mathbb{T},
L^1(\mathcal{M}) )$ be the Hardy space of $L^1 (\mathcal{M})$-valued
functions on $\mathbb{T}.$ If $f \in \mathrm{H}^1_c (\mathbb{T},
\mathcal{M} ),$ then
\begin{equation}\label{3.2}\| f \|_{\mathrm{H}^1} \leq \| f
\|_{\mathrm{H}^1_c}. \end{equation}Consequently, $\mathrm{H}^1_c
(\mathbb{T}, \mathcal{M} ) \subset \mathrm{H}^1 (\mathbb{T},
L^1(\mathcal{M}) )$ and is a Banach space under the norm $\| \cdot
\|_{\mathrm{H}^1_c}.$
\end{proposition}

\begin{proof}
As noted above, $\mathcal{M}_c$-atoms are all $L^1
(\mathcal{M})$-valued $2$-atoms. This implies the required result.
\end{proof}

\begin{definition}\label{df3.2}
We define\begin{equation*}\mathcal{H}^1_c (\mathbb{T}, \mathcal{M} )
= \left \{ f \in \mathrm{H}^1_c (\mathbb{T}, \mathcal{M} ): P(f) \in
\mathcal{H} (\mathbb{D}, L^1 ( \mathcal{M} ) ) \right
\},\end{equation*}equipped with the norm $\| f \|_{\mathrm{H}^1_c}.$
\end{definition}

By \eqref{3.2} one concludes that $\mathcal{H}^1_c (\mathbb{T},
\mathcal{M} ) \subset \mathcal{H}^1 (\mathbb{T}, L^1 ( \mathcal{M} )
),$ the Hardy space of $L^1 ( \mathcal{M} )$-valued functions on
$\mathbb{T},$ the Poisson integral of which are analytic in
$\mathbb{D}.$ Hence, $\mathcal{H}^1_c (\mathbb{T}, \mathcal{M} )$ is
a Banach space.

Similarly, we can define $\mathrm{H}^1_r (\mathbb{T}, \mathcal{M} )$
and $\mathcal{H}^1_r (\mathbb{T}, \mathcal{M} ).$

\begin{definition}\label{df3.3}
The Hardy space of operator-valued
analytic functions in the unit disc $\mathbb{D}$ is defined by
\begin{equation*}\mathcal{H}^1_{cr} (\mathbb{T}, \mathcal{M} )
= \mathcal{H}^1_c (\mathbb{T}, \mathcal{M} ) + \mathcal{H}^1_r
(\mathbb{T}, \mathcal{M} ), \end{equation*}equipped with the sum
norm\begin{equation*}\| f \|_{\mathcal{H}^1_{cr}} = \inf \left \{ \|
g \|_{\mathrm{H}^1_c} + \| h \|_{\mathrm{H}^1_r} \right
\},\end{equation*}where the infimum is taken over all $g \in
\mathcal{H}^1_c (\mathbb{T}, \mathcal{M} )$ and $h \in
\mathcal{H}^1_r (\mathbb{T}, \mathcal{M} )$ such that $f = g +h.$
\end{definition}

Clearly, $\mathcal{H}^1_{cr} (\mathbb{T}, \mathcal{M} )$ is a Banach
space. Let $\mathcal{H}^1 (\mathbb{D})$ denote the analytic Hardy
space in the unit disc $\mathbb{D}.$

\begin{proposition}\label{prop3.2}
If $f \in \mathcal{H}^1_{cr} (\mathbb{T}, \mathcal{M} )$ and $m \in
L^{\infty} (\mathcal{M}),$ then $\tau (m f) \in \mathcal{H}^1
(\mathbb{D})$ and\begin{equation}\label{3.3} \| \tau (m f)
\|_{\mathrm{H}^1} \leq \|m \|_{\infty} \| f \|_{\mathrm{H}^1_{cr}}.
\end{equation}
\end{proposition}

\begin{proof}
It is sufficient to prove that$$\| \tau (m a) \|_{\mathrm{H}^1} \leq
\|m \|_{\infty}$$for each $\mathcal{M}_c$-atom $a.$ Let $a$ be an
$\mathcal{M}_c$-atom supported in $I.$ By Lemma 2.2 we have$$\left (
\int_I |\tau (m a) |^2 d m \right )^{1/2} \leq \|m\|_{\infty} \tau
\left ( \int_I | a|^2 d m \right )^{1/2} \leq \|m\|_{\infty}
|I|^{-1/2}.$$This shows that $\tau (m a)/\| m \|_{\infty}$ is an
$2$-atom supported in $I.$ The proof is complete.
\end{proof}

\begin{remark}\label{re3.1}
The $\mathcal{M}_c$-atom is a noncommutative
analogue of the classical $2$-atom for $H^1.$ However, the
noncommutative analogues of classical $p$-atoms seem to be of no
meaningfulness.
\end{remark}

\section{Operator-valued BMOA space}

Let $f \in L^{\infty} (\mathcal{M}, L^2_c ( \mathbb{T}) ).$ Then, by
\eqref{2.6} one concludes that $\langle f, g \rangle \in
\mathcal{M}$ for each $g \in L^{\infty} ( \mathbb{T}).$ Thus, for a
subarc $I$ of $\mathbb{T}$ we can define the mean value $f_I$ of $f$
over $I$ by $f_I = \int_I f d m / |I| \in \mathcal{M}.$
Set\begin{equation*}\| f \|_{*, c} = \sup_I \left \| \left (
\frac{1}{|I|} \int_I | f - f_I |^2 d m \right )^{1/2} \right
\|_{\mathcal{M}},\end{equation*}where the supremum is taken over all
subarcs $I$ of $\mathbb{T}.$ We
define\begin{equation*}\mathrm{BMO}_c (\mathbb{T}, \mathcal{M}) =
\{f \in L^{\infty} (\mathcal{M}, L^2_c ( \mathbb{T}) ): \| f \|_{*,
c} < \infty \}\end{equation*}equipped with the
norm\begin{equation*}\| f \|_{\mathrm{BMO}_c} = \left \|
\int_{\mathbb{T}} f d m \right \|_{\mathcal{M}} + \| f \|_{*,
c}.\end{equation*}

\begin{proposition}\label{prop4.1}
If $f \in L^{\infty} (\mathcal{M}, L^2_c ( \mathbb{T}) ),$
then\begin{equation}\label{4.1} \| f \|_{*, c} \leq \sup_I \left (
\frac{1}{|I|} \int_I \| f - f_I \|^2_{\mathcal{M}} d m \right
)^{1/2}, \end{equation}and
\begin{equation}\label{4.2}\Big \| \int_{\mathbb{T}} f d m \Big
\|_{\mathcal{M}} \leq \| f \|_{L^{\infty}_c} \leq \Big \|
\int_{\mathbb{T}} f d m \Big \|_{\mathcal{M}} + \| f \|_{*, c}.
\end{equation}

Consequently, $\mathrm{BMO}(\mathbb{T}, \mathcal{M}) \subset
\mathrm{BMO}_c (\mathbb{T}, \mathcal{M})$ and $\mathrm{BMO}_c
(\mathbb{T}, \mathcal{M})$ is a Banach space. Here,
$\mathrm{BMO}(\mathbb{T}, \mathcal{M})$ is the $\mathrm{BMO}$ space
of $\mathcal{M}$-valued functions on $\mathbb{T}.$
\end{proposition}

\begin{proof}
Let $\mathbb{H}$ be the Hilbert space on which
$\mathcal{M}$ acts. For every $h \in \mathbb{H}$ one
has\begin{equation*}\int_I \| f(h) - f(h)_I \|^2 d m = \int_I \Big
\langle | f - f_I |^2 (h), h \Big \rangle d m = \Big \langle \Big (
\int_I | f - f_I |^2 d m \Big ) (h), h \Big \rangle,
\end{equation*}which yields that\begin{equation}\label{4.3} \| f
\|_{*, c} = \sup_{h \in \mathbb{H}, \|h \| \leq 1} \left \| f (h)
\right \|_{\mathrm{BMO} (\mathbb{T}, \mathbb{H})},
\end{equation}where $\mathrm{BMO} (\mathbb{T}, \mathbb{H})$
is the $\mathbb{H}$-valued BMO space on $\mathbb{T}.$ This concludes
\eqref{4.1} and so, $\mathrm{BMO}(\mathbb{T}, \mathcal{M}) \subset
\mathrm{BMO}_c (\mathbb{T}, \mathcal{M}).$

Since\begin{equation*}\| f \|^2_{L^{\infty}_c} = \Big \|
\int_{\mathbb{T}} | f |^2 d m \Big \|_{\mathcal{M}} = \sup_{\| h \|
\leq 1} \Big \langle \int_{\mathbb{T}} | f |^2 d m (h), h \Big
\rangle = \sup_{\| h \| \leq 1} \int_{\mathbb{T}} \| f (h) \|^2 d
m,\end{equation*}it is concluded that\begin{equation*}\| f
\|_{L^{\infty}_c} = \sup_{\| h \| \leq 1} \Big ( \int_{\mathbb{T}}
\| f (h) \|^2 d m \Big )^{1/2} \geq \sup_{\| h \| \leq 1}
\int_{\mathbb{T}} \| f (h) \| d m \geq \Big \| \int_{\mathbb{T}} f d
m \Big \|_{\mathcal{M}},
\end{equation*}and\begin{equation*}\label{}\begin{split}\| f
\|_{L^{\infty}_c} &=  \sup_{\| h \| \leq 1} \Big ( \int_{\mathbb{T}}
\| f (h) \|^2 d m \Big )^{1/2}\\
& \leq \sup_{\| h \| \leq 1} \Big ( \int_{\mathbb{T}} \| f (h) -
\int_{\mathbb{T}} f(h) d m \|^2 d m \Big )^{1/2} + \sup_{\| h \|
\leq 1} \Big \| \int_{\mathbb{T}} f(h) d m \Big \|\\
& \leq \Big \| \int_{\mathbb{T}} f d m \Big \|_{\mathcal{M}} + \| f
\|_{*, c}. \end{split}\end{equation*}This completes the proof of
\eqref{4.2}.
\end{proof}

\begin{definition}\label{df4.1}
We define\begin{equation*}\mathrm{BMOA}_c (\mathbb{T}, \mathcal{M})
= \{ f \in \mathrm{BMO}_c (\mathbb{T}, \mathcal{M}): P[f] \in
\mathcal{H} (\mathbb{D},\mathcal{M} ) \},\end{equation*}equipped
with the norm $\| \cdot \|_{\mathrm{BMO}_c}.$
\end{definition}

By \eqref{4.1} and \eqref{4.3} we conclude that $\mathrm{BMOA}
(\mathbb{T}, \mathcal{M}) \subset \mathrm{BMOA}_c (\mathbb{T},
\mathcal{M})$ and $\mathrm{BMOA}_c (\mathbb{T}, \mathcal{M})$ is
Banach space.

For $f \in \mathrm{BMOA}_c (\mathbb{T}, \mathcal{M})$ we introduce
the set of functions$$M_f = \{ g \in \mathcal{H}
(\mathbb{D},\mathcal{M} ): g = f \circ \psi - f \circ \psi (0), \psi
\in \mathrm{Aut} (\mathbb{D}) \},$$and for $0 < \gamma < 1,$ we
denote by $f_{\gamma}$ the dilated function defined for $|z| < 1/
\gamma$ by $f_{\gamma} (z) = f( \gamma z).$ Also, for $f \in
L^{\infty} (\mathcal{M}, L^2_c ( \mathbb{T}) )$ set$$\| f \|_{**, c}
= \sup_{z \in \mathbb{D}} \left \| \left ( \int_{\mathbb{T}} | f - f
(z) |^2 d m_z \right )^{1/2} \right \|_{\mathcal{M}},$$where $ d m_z
(t) = P_z (t ) d m (t)$ for $z \in \mathbb{D}.$

\begin{proposition}\label{prop4.2}
If $f \in \mathrm{BMOA}_c (\mathbb{T}, \mathcal{M}),$
then\begin{equation}\label{4.4} \| f \|_{*, c} \asymp \sup_{g \in
M_f} \sup_{0 < \gamma < 1} \left \| g_{\gamma} \right
\|_{L^{\infty}_c} \asymp \| f \|_{**, c}. \end{equation}
\end{proposition}

\begin{proof}
Let $\mathbb{H}$ be the Hilbert space on which
$\mathcal{M}$ acts. By merely reproducing the proof in the case of
scalars (for details, see \cite{BII}) we can obtain \eqref{4.4} for
$\mathbb{H}$- valued functions. Then, by \eqref{4.3} we
have\begin{equation*}\label{}\begin{split}\| f \|_{*, c} & = \sup_{h
\in \mathbb{H}, \|h \| \leq 1} \left \| f
(h) \right \|_{\mathrm{BMO} (\mathbb{T}, \mathbb{H})} \\
& \asymp \sup_{ \|h \| \leq 1} \sup_{g \in M_f} \sup_{0< \gamma <1}
\left ( \int_{\mathbb{T}} \| g_{\gamma} (h) \|^2 d m \right
)^{1/2} \\
& = \sup_{g \in M_f} \sup_{0 < \gamma <1} \left \| g_{\gamma} \right
\|_{L^{\infty}_c}.\end{split}
\end{equation*}This concludes the first equivalence.

Similarly, for $f \in \mathrm{BMOA}_c (\mathbb{T}, \mathcal{M})$ we
have\begin{equation*}\label{}\begin{split}\| f \|_{*, c} & = \sup_{h
\in \mathbb{H}, \|h \| \leq 1} \left \| f(h) \right \|_{\mathrm{BMO} (\mathbb{T}, \mathbb{H})} \\
& \asymp \sup_{ \|h \| \leq 1} \sup_{z \in \mathbb{D}} \left (
\int_{\mathbb{T}} \| f(h) - f(h) (z) \|^2 d m_z \right )^{1/2}\\
& = \sup_{z \in \mathbb{D}} \left \| \left ( \int_{\mathbb{T}} | f -
f (z) |^2 d m_z \right )^{1/2} \right \|_{\mathcal{M}}.
\end{split}\end{equation*}This proves that $\| f \|_{*, c} \asymp \| f \|_{**,
c}$ and completes the proof.
\end{proof}

Similarly, we can define $\mathrm{BMO}_r (\mathbb{T}, \mathcal{M})$
and $\mathrm{BMOA}_r (\mathbb{T}, \mathcal{M}),$ by letting $$\| f
\|_{*, r} = \| f^* \|_{*, c}\; \text{and}\; \| f \|_{\mathrm{BMO}_r}
= \| f^* \|_{\mathrm{BMO}_c},$$ respectively. Evidently, we have

\begin{proposition}\label{prop4.3}
If $f \in \mathrm{BMOA}_r (\mathbb{T},
\mathcal{M}),$ then\begin{equation}\label{4.5}
 \| f \|_{*, r} \asymp \sup_{g \in M_f}
\sup_{0 < \gamma <1} \left \| g_{\gamma} \right \|_{L^{\infty}_r}
\asymp \| f \|_{**, r}.
\end{equation}Here, $\| f \|_{**, r} = \| f^* \|_{**, c}.$
\end{proposition}

We define\begin{equation*}\mathrm{BMO}_{cr} (\mathbb{T},
\mathcal{M}) = \mathrm{BMO}_c (\mathbb{T}, \mathcal{M}) \cap
\mathrm{BMO}_r (\mathbb{T}, \mathcal{M}),\end{equation*}quipped with
the norm $\| f \|_{\mathrm{BMO}_{cr}}= \max \left \{ \| f
\|_{\mathrm{BMO}_c}, \| f \|_{\mathrm{BMO}_r} \right \}.$ As shown
above, $\mathrm{BMO}(\mathbb{T}, \mathcal{M}) \subset
\mathrm{BMO}_{cr} (\mathbb{T}, \mathcal{M}).$

\begin{definition}\label{df4.2}
The operator-valued $\mathrm{BMOA}$ space on $\mathbb{T}$ is defined
as follows:\begin{equation*}\mathrm{BMOA}_{cr} (\mathbb{T},
\mathcal{M}) = \mathrm{BMOA}_c (\mathbb{T}, \mathcal{M}) \cap
\mathrm{BMOA}_r (\mathbb{T}, \mathcal{M}),\end{equation*}quipped
with the norm $\| f \|_{\mathrm{BMO}_{cr}}.$
\end{definition}

$\mathrm{BMOA}_{cr} (\mathbb{T}, \mathcal{M})$ is evidently a Banach
space under the norm $\| \cdot \|_{\mathrm{BMO}_{cr}}.$ Moreover,
$\mathrm{BMOA}(\mathbb{T}, \mathcal{M}) \subset \mathrm{BMOA}_{cr}
(\mathbb{T}, \mathcal{M}).$

It is well-known that BMO spaces are related to the so-called
Carleson measures (e.g., see \cite{G, S93}). In the sequel, we do
this in the operator-valued setting on the circle, with an emphasis
on the `analytical' aspect.

For $t \in \mathbb{T}$ and $\delta > 0$ we introduce the
set$$\hat{I} (t_0, \delta) = \{rt \in \mathbb{D}: 1- \delta \leq r <
1, |t - t_0| < \delta \},$$whose closure intersects $\mathbb{T}$ at
the subarc $I(t_0, \delta) = \{t \in \mathbb{T}: |t-t_0|< \delta
\}.$ $\hat{I} (t, \delta)$ is said to be a Carleson tube at $t.$

\begin{definition}\label{df4.3}
An $\mathcal{M}_+$-valued Borel measure $\nu$ in $\mathbb{D}$ is
called a Carleson measure if there exists a constant $C>0$ such
that$$\| \nu (\hat{I} (t, \delta)) \|_{\mathcal{M}} \leq C
\delta,$$for all $t \in \mathbb{T}$ and $\delta > 0.$
Set\begin{equation*}\| \nu \|_c = \sup_I \frac{\| \nu (\hat{I} (t,
\delta)) \|_{\mathcal{M}}}{\delta},\end{equation*}which is called
the Carleson norm of $\nu.$
\end{definition}

Since $|I (t, \delta)| \asymp \delta$ ($0 < \delta < 1$) and
$\hat{I} (t, \delta) = \mathbb{D}$ provided $\delta > 1,$ it is
concluded that an $\mathcal{M}_+$-valued Borel measure $\nu$ in
$\mathbb{D}$ is a Carleson measure if and only
if\begin{equation}\label{4.6} \sup \Big \{ \frac{\| \nu (\hat{I} (t,
\delta)) \|_{\mathcal{M}}}{|I (t, \delta)|}: t \in \mathbb{T},
\delta > 0 \Big \} < \infty. \end{equation}

\begin{proposition}\label{prop4.4}
An $\mathcal{M}_+$-valued Borel measure $\nu$ in $\mathbb{D}$ is a
Carleson measure if and only if$$\mathcal{N}(\nu) = \sup_{z \in
\mathbb{D}} \left \| \int_{\mathbb{D}} P_z(w) d \nu (w) \right
\|_{\mathcal{M}} < \infty.$$The constant $\mathcal{N}(\nu)$
satisfies\begin{equation*}C_1 \mathcal{N}(\nu) \leq \| \nu \|_c \leq
C_2 \mathcal{N}(\nu) \end{equation*}for absolute constants $C_1$ and
$C_2.$
\end{proposition}

The proof is the same as in the scalar case (e.g., see Lemma 3.3 in
\cite{G}) and omitted. Proposition \ref{prop4.4} shows the
conformally invariant character of Carleson measures.

The following theorem is the operator-valued version of one of the
most fundamental results in the theory of BMO spaces, which
characterizes BMO spaces in terms of Carleson measures.

\begin{theorem}\label{th4.1}
Let $f \in L^{\infty} (\mathcal{M}, L^2_c (
\mathbb{T}) ).$ Then, the following assertions are equivalent:
\begin{enumerate}[{\rm (1)}]

\item $f$ is in $\mathrm{BMO}_c (\mathbb{T}, \mathcal{M}).$

\item $ d \nu_f (z) = | \nabla f(z) |^2 (1-|z|^2 ) d x d y$ is a Carleson
measure in $\mathbb{D}.$

\item $ d \lambda_f (z) = | \nabla f(z) |^2 \log \frac{1}{|z|} d x d y$ is a Carleson
measure in $\mathbb{D}.$
\end{enumerate}
In this case, \begin{equation*}\| \nu_f \|_c \asymp \| f \|^2_{*,c}
\asymp \| \lambda_f \|_c.\end{equation*}

Consequently, if $f \in \mathcal{H}^{\infty} (\mathcal{M}, L^2_{cr}
( \mathbb{T}) )$ then $f \in \mathrm{BMOA}_{cr} (\mathbb{T},
\mathcal{M})$ if and only if both $(1 - |z|^2 )| f'(z) |^2 d x d y$
and $(1 - |z|^2 )| f'(z)^* |^2 d x d y$ are $\mathcal{M}_{+}$-valued
Carleson measures in $\mathbb{D}$ if and only if $$| f'(z) |^2 \log
\frac{1}{|z|} d x d y\; \text{and}\; | f'(z)^* |^2 \log
\frac{1}{|z|} d x d y$$ are both $\mathcal{M}_{+}$-valued Carleson
measures in $\mathbb{D}.$
\end{theorem}

\begin{proof}
By Lemma \ref{le2.3}, Propositions \ref{prop4.2} and \ref{prop4.4}
we get the equivalence of (1) and (2) and $\| \nu_f \|_c \asymp \| f
\|^2_{*,c}.$ What remains to prove is the equivalence of (2) and
(3). Half of this task is trivial because the inequality $1-|z|^2
\leq 2 \log (1/|z|)$ shows that $\| \nu_f \|_c \leq 2 \| \lambda_f
\|_c.$ For $|z| > 1/4$ we have the reverse inequality $\log (1/ |z|)
\leq C (1-|z|^2),$ which shows that $\| \lambda_f (\hat{I}) \| \leq
C \| \nu_f (\hat{I}) \| $ for subarcs $I = I (t, \delta)$ provided
$\delta \leq 3/4,$ because $|z| > 1 - \delta \geq 1/4$ for all $z
\in \hat{I}.$ Then, to prove that the equivalence of (2) and (3) it
suffices to prove $\| \lambda_f ( |z| \leq 1/4) \| \leq C \| \nu_f (
|z| \leq 1/2) \|.$ However, as shown in the proof of Lemma
\ref{le2.3} we have\begin{equation*}|\nabla f (z) |^2 \leq C
\int_{|z| < 1/2}|\nabla f (w) |^2 (1-|w|^2) dx dy = C \nu_f (|z|
\leq 1/2)
\end{equation*}for all $|z| < 1/4.$ Hence,
\begin{equation*}\lambda_f  ( |z| \leq 1/4)= \int_{|z| \leq 1/4}
|\nabla f (z) |^2 \log \frac{1}{|z|} dx dy \leq C \nu_f (|z| \leq
1/2).\end{equation*}This gives that $\| \lambda_f  ( |z| \leq 1/4)
\| \leq C \| \nu_f ( |z| \leq 1/2) \|$ and therefore $\| \lambda_f
\|_c \leq C \| \nu_f \|_c.$
\end{proof}

\section{Operator-valued $\mathrm{H}^1$-BMOA duality}

In this section we show the operator-valued $\mathrm{H}^1$-BMOA
duality.

\begin{theorem}\label{th5.1}
We have \begin{equation*}\mathrm{H}^1_c (\mathbb{T}, \mathcal{M} )^*
= \mathrm{BMO}_c (\mathbb{T}, \mathcal{M} ),~~\mathcal{H}^1_c
(\mathbb{T}, \mathcal{M} )^* = \mathrm{BMOA}_c (\mathbb{T},
\mathcal{M} ).\end{equation*}Similarly, $\mathrm{H}^1_r (\mathbb{T},
\mathcal{M} )^* = \mathrm{BMO}_r (\mathbb{T}, \mathcal{M} )$ and
$\mathcal{H}^1_r (\mathbb{T}, \mathcal{M} )^* = \mathrm{BMOA}_r
(\mathbb{T}, \mathcal{M} ).$ Consequently, $\mathcal{H}^1_{cr}
(\mathbb{T}, \mathcal{M} )^* = \mathrm{BMOA}_{cr} (\mathbb{T},
\mathcal{M} ).$
\end{theorem}

\begin{proof}
Let $g \in \mathrm{BMO}_c (\mathbb{T}, \mathcal{M} ).$ For any
$\mathcal{M}_c$-atom $a$ supported in $I,$ we have by Lemma
\ref{le2.1}
\begin{equation*}\label{}\begin{split}\Big | \tau \Big (
\int_{\mathbb{T}} g^* a d m \Big ) \Big | & = \Big | \tau \Big (
\int_I (g -
g_I )^* a d m \Big ) \Big | \\
& \leq \| a \|_{L^1 (\mathcal{M}, L^2_c (
I) )} \| g - g_I \|_{L^{\infty} (\mathcal{M}, L^2_c ( I) )} \\
& \leq \| g \|_{*, c}.\end{split}\end{equation*}On the other hand,
for any $a \in L^1 (\mathcal{M})$ with $\| a \| \leq 1$ we
have$$\Big | \tau \Big ( \int_{\mathbb{T}} g^* a d m \Big ) \Big |
\leq \Big \| \int_{\mathbb{T}} g^* d m \Big \| \tau (|a|) \leq \Big
\| \int_{\mathbb{T}} g d m \Big \|.$$ Thus, we deduce
that\begin{equation}\label{5.1}\Big | \tau \Big ( \int_{\mathbb{T}}
g^* f d m \Big )\Big | \leq \| g \|_{\mathrm{BMO}_c} \|
f\|_{\mathrm{H}^1_c}
\end{equation}for all $f \in \mathrm{H}^1_c (\mathbb{T},
\mathcal{M} ).$ Hence, $\mathrm{BMO}_c (\mathbb{T}, \mathcal{M} )
\subset \mathrm{H}^1_c (\mathbb{T}, \mathcal{M} )^*.$

Conversely, let $l \in \mathrm{H}^1_c (\mathbb{T}, \mathcal{M} )^*.$
Since $L^1 (\mathcal{M}, L^2_c ( \mathbb{T}) ) \subset
\mathrm{H}^1_c (\mathbb{T}, \mathcal{M} ),$ $l$ induces a continuous
functional on $L^1 (\mathcal{M}, L^2_c ( \mathbb{T} ) )$ with the
norm smaller than $ \| l \|_{( \mathrm{H}^1_c)^*},$ that is, $| l (f
) | \leq \| l \|_{( \mathrm{H}^1_c)^*} \| f \|_{L^1_c}$ for all $f
\in L^1 (\mathcal{M}, L^2_c ( \mathbb{T} ) ).$ Hence, by the
Hahn-Banach extension theorem there exists a unique $g \in
L^{\infty} (\mathcal{M}, L^2_c ( \mathbb{T} ) )$ with $\| g
\|_{L^{\infty}_c} \leq  \| l \|_{( \mathcal{H}^1_c)^*}$ such
that\begin{equation}\label{5.2}l (f ) = \tau \Big (
\int_{\mathbb{T}} g^* f d m \Big )
\end{equation}for all $f \in L^1 (\mathcal{M}, L^2_c (
\mathbb{T} ) ).$ We need to show that $g \in \mathrm{BMO}_c
(\mathbb{T}, \mathcal{M} ).$

Let $\mathbb{H}$ be the Hilbert space on which $\mathcal{M}$ acts.
Recall that the predual space $L^1 ( \mathcal{M} )$ of $\mathcal{M}$
is the quotient space $L^1 ( \mathcal{B} (\mathbb{H})
)/\mathcal{M}_{\perp},$ where $\mathcal{M}_{\perp} = \{ a \in L^1 (
\mathcal{B} (\mathbb{H}) ): \mathrm{Tr} (a b) = 0, \forall b \in
\mathcal{M} \}$ is the pre-annihilator of $\mathcal{M}.$ The
quotient map is dented by $\pi: L^1 ( \mathcal{B} (\mathbb{H}) )
\rightarrow L^1 ( \mathcal{B} (\mathbb{H}) ) / \mathcal{M}_{\perp}.$
For $u, v \in \mathbb{H}$ we define $|u \rangle \langle v |$ by $|u
\rangle \langle v | (h ) = \langle h, v \rangle u$ for all $h \in
\mathbb{H}.$ Then, for every $a \in \mathcal{M}$ one has$$\tau (a
\pi [ |u \rangle \langle v | ] ) = \tau ( \pi [ a |u \rangle \langle
v | ] ) = \mathrm{Tr} (a |u \rangle \langle v |) = \langle a (u ), v
\rangle.$$

Consequently, by \eqref{4.3} and the Hilbert space-valued
$\mathrm{H}^1$-BMOA duality theorem we
have\begin{equation*}\begin{split}\| g \|_{*, c}
& = \sup_{h \in \mathbb{H}, \| h \| = 1} \| g (h) \|_{\mathrm{BMO}(\mathbb{T}, \mathbb{H})}\\
& = \sup_{\| h \| = 1} \sup_{f \in \mathrm{H}^1 (\mathbb{T},
\mathbb{H}), \| f \|_{\mathrm{H}^1} \leq 1} \Big |
\int_{\mathbb{T}} \langle g(h), f \rangle d m  \Big |\\
& = \sup_{\| h \| = 1} \sup_{\| f \|_{\mathrm{H}^1} \leq 1} \Big |
\int_{\mathbb{T}} \mathrm{Tr}[g | h \rangle \langle f | ] d m \Big |\\
& = \sup_{\| h \| = 1} \sup_{\| f \|_{\mathrm{H}^1} \leq 1 } \Big |
\tau \Big ( \int_{\mathbb{T}} g \pi [ | h
\rangle \langle f | ]  d m \Big )  \Big |\\
& \leq \sup_{\| f \|_{\mathrm{H}^1_c} \leq 1 } \Big | \tau \Big (
\int_{\mathbb{T}} g f^* d m \Big ) \Big |\\
& = \| l \|_{( \mathrm{H}^1_c )^*}.
\end{split}\end{equation*}
Also,\begin{equation*}\begin{split}\Big \| \int_{\mathbb{T}} g d m
\Big \| & = \sup_{y \in L^1 (\mathcal{M}), \|y \|_1 \leq 1} \Big |
\tau \Big ( \int_{\mathbb{T}} g^* y d m \Big ) \Big |\\
& \leq \sup_{f \in L^1_c, \|f \|_{L^1_c} \leq 1} \Big | \tau \Big (
\int_{\mathbb{T}} g^* f d m \Big ) \Big | \leq \| l \|_{(
\mathrm{H}^1_c )^*}.\end{split}\end{equation*}This shows that $\| g
\|_{\mathrm{BMO}_c} \leq 2 \| l \|_{( \mathrm{H}^1_c )^*}$ and
concludes the first assertion.

By \eqref{5.1}, each $g \in \mathrm{BMOA}_c (\mathbb{T}, \mathcal{M}
)$ induces evidently a bounded linear functional on $\mathcal{H}^1_c
(\mathbb{T}, \mathcal{M} )$ via \eqref{5.2}. Conversely, if $l \in
\mathcal{H}^1_c (\mathbb{T}, \mathcal{M} )^*,$ then by the
Hahn-Banach extension theorem, $l$ can be extended to a bounded
linear functional on $\mathrm{H}^1_c (\mathbb{T}, \mathcal{M} )$ and
hence there exists a $g \in L^{\infty} (\mathcal{M}, L^2_c (
\mathbb{T} ) )$ such that \eqref{5.2} holds true for all (analytic)
polynomials $f= \sum^n_{k=1} b_k z^k,$ where $b_k \in L^1
(\mathcal{M}).$ This concludes that\begin{equation}\label{5.3}l (P )
= \tau \Big ( \int_{\mathbb{T}} \mathfrak{C}[g]^* P d m \Big )
\end{equation}for all polynomials $P (z) = \sum^n_{k=1} b_k
z^k,$ where $b_k \in L^1 (\mathcal{M}).$ By merely reproducing the
above proof we can prove that $\mathfrak{C}[g] \in \mathrm{BMOA}_c
(\mathbb{T}, \mathcal{M} ).$
\end{proof}

\section{Area integral characterizations}

We define for each $f \in L^1 (\mathcal{M}, L^2_c ( \mathbb{T} ) )
:$

(i)  The {\it column} area function\begin{equation*}\label{} [A_c
(f, \alpha)] (t) = \left ( \int_{\Gamma_{\alpha} (t)} | \nabla f (z)
|^2 dx dy \right )^{1/2}, ~~(t \in \mathbb{T}),\end{equation*}

(ii)  The {\it column} Littlewood-Paley $g$-function
\begin{equation*}\label{}[g_c
(f)] (t)  = \left ( \int^1_0 | \nabla f (r t) |^2 (1-r^2) d r \right
)^{1/2},~~(t \in \mathbb{T}).\end{equation*}

Similarly, we define the {\it row} area function $A_r (f, \alpha) =
A_c (f^*, \alpha)$ and {\it row} Littlewood-Paley $g$-function $g_r
(f) = g_c (f^*).$

Also, we need two technical variants of $A_c (f, \alpha)$ and $g_c
(f)$ as following:\begin{equation*}\label{}[A_c (f, \alpha)] (t,
\delta) = \Big ( \int_{\Gamma_{\alpha} (t, \delta)} | \nabla f (z)
|^2 dx dy \Big )^{1/2},\end{equation*}where $\Gamma_{\alpha} (t,
\delta) = \{z \in \mathbb{D}: |t -z| < \alpha (1-|z|), |z|< \delta
\},$ and\begin{equation*}\label{}[g_c (f)] (t, \delta) = \Big (
\int^{\delta}_0 | \nabla f (r t) |^2 (1-r^2) d r \Big
)^{1/2},\end{equation*}for $0 < \delta \leq 1, t \in \mathbb{T}.$

For simplicity, we usually denote by $A_c (f)= A_c (f, \alpha)$ in
the sequel.

\begin{lemma}\label{le6.1}
There is a constant $C>0$ such that\begin{equation*} \label{}[g_c
(f)] (t, \delta) \leq C [A_c (f)] (t, (1 + \delta)/2),~~ (0 < \delta
\leq 1,~ t \in \mathbb{T}),\end{equation*}for all $f \in L^1
(\mathcal{M}, L^2_c ( \mathbb{T} ) ).$ In particular, $g_c (f) \leq
C A_c (f).$
\end{lemma}

\begin{proof}
It suffices to prove the associated inequality
for the case of $t=1.$ Since $\Gamma_{\alpha} (1, (1+ \delta)/2)$
contains a small disc centered at $0,$ there exists a constant $0 <
c_{\alpha} < 1$ such that for each $0 < r <
\delta,$\begin{equation*}D_r \triangleq \{z \in \mathbb{D}: | z - r|
< c_{\alpha}[(1+\delta)/2 - r] \} \subset \Gamma_{\alpha} (1, (1+
\delta)/2).\end{equation*}The harmonicity of $\nabla f$ gives
that\begin{equation*}\nabla f (r) = \frac{1}{c^2_{\alpha}\pi
[(1+\delta)/2 - r]^2} \int_{D_r} \nabla f (z) dx d
y.\end{equation*}This follows from \eqref{2.8}
that\begin{equation*}|\nabla f (r)|^2 \leq \frac{C}{(1- r)^2}
\int_{D_r} | \nabla f (z)|^2 dx d
y.\end{equation*}Hence,\begin{equation*}[g_c (f)] (1, \delta) =
\int^{\delta}_0 |\nabla f (r)|^2 (1-r^2) dr \leq C \int^{\delta}_0
\frac{d r}{1- r} \int_{D_r} | \nabla f (z)|^2 dx d
y.\end{equation*}Since\begin{equation*}(1-|z|)/(1 + c_{\alpha}) < 1-
r < (1- |z|)/(1 - c_{\alpha})\end{equation*}for $z \in D_r,$ by
Fubini's theorem we have\begin{equation*}[g_c (f)] (1, \delta) \leq
C \int_{\Gamma_{\alpha} \big ( 1, \frac{1+\delta}{2} \big )} |
\nabla f (z)|^2 \int^{\frac{1- |z|}{1 - c_{\alpha}}
}_{\frac{1-|z|}{1 + c_{\alpha}}} \frac{d r}{r} dx d y \leq C [A_c
(f)] \big ( 1, \frac{1 + \delta}{2} \big ).\end{equation*}This
completes the proof.
\end{proof}

\begin{theorem}\label{th6.1}
Let $f \in L^1 (\mathcal{M}, L^2_c (
\mathbb{T}) ).$ The following assertions are equivalent:
\begin{enumerate}[{\rm (1)}]

\item  $f \in \mathrm{H}^1_c (\mathbb{T}, \mathcal{M} ).$

\item  $A_c (f) \in L^1(L^{\infty}(\mathbb{T}) \otimes \mathcal{M}).$

\item  $g_c (f) \in L^1(L^{\infty}(\mathbb{T}) \otimes \mathcal{M}).$

\end{enumerate}
In this case,\begin{equation}\label{6.1}\|f\|_{\mathrm{H}^1_c}
\asymp \| f(0) \|_{L^1 (\mathcal{M})} + \| A_c (f) \|_{L^1} \asymp
\| f(0) \|_{L^1 (\mathcal{M})} + \| g_c (f) \|_{L^1}
\end{equation}for all $f \in \mathrm{H}^1_c (\mathbb{T},
\mathcal{M} ).$

Consequently, if $f \in \mathcal{H}^1 (\mathcal{M}, L^2_c (
\mathbb{T}) ),$ then $f \in \mathcal{H}^1_c (\mathbb{T}, \mathcal{M}
)$ if and only if $A_c (f) \in L^1(\mathbb{T}, L^1(\mathcal{M}))$ if
and only if $g_c (f) \in L^1(\mathbb{T}, L^1(\mathcal{M})).$

The same statements hold also true for $\mathrm{H}^1_r (\mathbb{T},
\mathcal{M} )$ and $\mathcal{H}^1_r (\mathbb{T}, \mathcal{M} ).$
\end{theorem}

\begin{proof}
We will show that $(1) \Rightarrow (2),$ $(2)
\Rightarrow (3),$ and finally $(3) \Rightarrow (1).$

$(1) \Rightarrow (2).$~~We will show that there exists $C>0$ such
that\begin{equation}\label{6.2}\| f (0) \|_{L^1 (\mathcal{M})} + \|
A_c (f) \|_{L^1} \leq C \|f\|_{\mathrm{H}^1_c}\end{equation}for all
$f \in \mathrm{H}^1_c (\mathbb{T}, \mathcal{M} ).$ Since $$\|
 f (0) \|_{L^1 (\mathcal{M})} \leq C \| f
\|_{\mathrm{H}^1} \leq  C \| f \|_{\mathrm{H}^1_c}$$ by Proposition
\ref{prop3.1}, it suffices to show that there exists $C>0$ such that
$$\| A_c (a) \|_{L^1} \leq C$$ for all $\mathcal{M}_c$-atom $a.$

Given an $\mathcal{M}_c$-atom $a$ supported in $I = I (t_0,
\delta).$ By Lemma \ref{le2.2} and Lemma \ref{le2.4} we
have\begin{equation*}\label{} \tau \Big ( \int_{2 I}  A_c (a) d m
\Big ) \leq |2 I|^{1/2} \Big ( \int_{2 I} \big [ \tau ( | A_c ( a) |) \big ]^2 d m
\Big )^{1/2} \leq C |I|^{1/2} \| a \|_{L^1_c} \leq
C,\end{equation*}where $2 I = I (t_0, 2 \delta ),$ and by
\eqref{2.8}\begin{equation*}\label{}\begin{split} \tau \Big (&
\int_{\mathbb{T}\setminus 2 I} A_c ( a ) d m \Big )\\
& = \int_{\mathbb{T}\setminus 2 I} \tau \Big [ \Big (
\int_{\Gamma_{\alpha}
(t)} | \nabla a (z) |^2 dx dy \Big )^{\frac{1}{2}} \Big ] d m (t) \\
& =  \int_{\mathbb{T}\setminus  2 I} \tau \Big [ \Big (
\int_{\Gamma_{\alpha} (t)} \Big | \int_I \big ( \nabla[P_z(s ) -
P_z( t_0 ) ] \big ) a (s) d m (s) \Big |^2 dx dy \Big )^{\frac{1}{2}} \Big ] d m (t)\\
& \leq \| a \|_{L^1_c} \int_{\mathbb{T}\setminus  2 I} \Big (
\int_{\Gamma_{\alpha} (t)} \int_I \Big |  \nabla[P_z(s) - P_z( t_0 )
] \Big |^2 d m(s)  dx
dy \Big )^{\frac{1}{2}} d m (t) \\
& \leq \int_{\mathbb{T}\setminus 2 I} \Big ( \int_{\Gamma_{\alpha}
(t)} \sup_{s \in I} \big | \nabla [ P_z( s ) - P_z(t_0 )] \big |^2
dx dy \Big )^{\frac{1}{2}} d m (t).\end{split}\end{equation*} An
immediate computation yields
that\begin{equation*}\begin{split}\frac{\partial P_z (t)}{\partial
x} & = - \frac{2x}{|1- \bar{t} z|^2} +
\frac{2(1-|z|^2)\mathrm{Re}(t- z)}{|1-\bar{t}z|^4},\\
\frac{\partial P_z (t)}{\partial y} & =- \frac{2y}{|1- \bar{t} z|^2}
+ \frac{2(1-|z|^2)\mathrm{Im}( t- z)}{|1-\bar{t}z|^4},
\end{split}\end{equation*}and hence,\begin{equation*}\begin{split}
\left | \frac{\partial}{\partial \gamma} \Big [\frac{\partial P_z
(t_0 + \gamma t)}{\partial x} \Big ] \right |, \left |
\frac{\partial}{\partial \gamma} \Big [\frac{\partial P_z (t_0 +
\gamma t)}{\partial y} \Big ] \right | \leq \frac{C|t|}{| 1 - (t_0 +
\gamma t) \bar{z} |^3},
\end{split}\end{equation*}for $0 < \gamma <1.$ Then,\begin{equation*}\begin{split}
\big | \nabla [ P_z( s ) - P_z(t_0 )] \big |^2 & = \Big |
\frac{\partial P_z (s)}{\partial x} - \frac{\partial P_z
(t_0)}{\partial x} \Big |^2 +  \Big | \frac{\partial P_z
(s)}{\partial y} - \frac{\partial P_z (t_0)}{\partial y} \Big |^2\\
& \leq \frac{C |s - t_0|^4}{| 1 - [t_0 + \gamma (s- t_0)] \bar{z}
|^6 }, \end{split}\end{equation*}for some $0 < \gamma < 1.$
Since\begin{equation*}\label{}|s' - t_0| = \gamma |s - t_0 | \leq
\delta,\end{equation*}where $s'= t_0 + \gamma (s- t_0),$ it is
concluded that
\begin{equation*}|t-t_0| \leq 2 |t - s' | \leq C(|1 - t \bar{z}| + |1 - s' \bar{z}|) \leq
C(1-|z| + |1 -s' \bar{z}|) \leq C|1 - s' \bar{z}|\end{equation*}for
$t \in \mathbb{T} \setminus 2 I$ and $z \in \Gamma_{\alpha} (t),$
the first inequality is a consequence of the triangle inequality and
the hypotheses. Hence,\begin{equation*} \big | \nabla [ P_z( s ) -
P_z(t_0 )] \big |^2 \leq \frac{C \delta^4}{| 1 - t_0\bar{t} |^6 },
\end{equation*}for $s \in I(t_0, \delta), t \in \mathbb{T} \setminus 2 I$ and
$z \in \Gamma_{\alpha} (t).$ This concludes that
\begin{equation*} \tau \Big (
\int_{\mathbb{T} \setminus 2 I} A_c ( a ) d m \Big ) \leq
\int_{\mathbb{T} \setminus 2 I} \frac{C \delta^2}{| 1 - t_0 \bar{t}
|^3 } d m (t) \leq C,\end{equation*}which completes the proof of
that $(1) \Rightarrow (2).$

$(2) \Rightarrow (3).$~~This follows from Lemma \ref{le6.1}.

$(3) \Rightarrow (1).$~~Since $\| f(0) \|_{\mathrm{H}^1_c} \leq \|
f(0) \|_{L^1 (\mathcal{M})},$ by Theorem \ref{th5.1} it suffices to
show that\begin{equation*}\Big | \tau \Big ( \int_{\mathbb{T}} f g^*
d m \Big ) \Big | \leq C \left \| g_c (f) \right \|_{L^1(\mathbb{T},
L^1(\mathcal{M}))} \left \| g \right \|_{\mathrm{BMO_c}
},\end{equation*}for all $f \in L^1_0 (\mathcal{M}, L^2_c (
\mathbb{T} ) )$ and $g \in L^{\infty}_0 (\mathcal{M}, L^2_c (
\mathbb{T} ) ).$ By Lemma \ref{le2.3} and the Cauchy-Schwarz
inequality, we have\begin{equation*}\begin{split}&\Big |  \tau \Big
( \int_{\mathbb{T}} f g^* d m \Big ) \Big |\\
& = \frac{1}{\pi} \left
| \tau \Big ( \int_{\mathbb{D}} \nabla f (z)
(\nabla g)^* (z) \log \frac{1}{|z|} dx dy  \Big ) \right | \\
& = \frac{2}{\pi} \left | \tau \Big ( \int^1_0 r \log \frac{1}{r} dr
\int_{\mathbb{T}} \nabla f (r t) (\nabla g)^* (r t) d m(t)  \Big ) \right |\\
& \leq \frac{2}{\pi} \left [ \tau \Big ( \int^1_0 r \log \frac{1}{r}
dr \int_{\mathbb{T}} [g_c (f )
(r, t)]^{- \frac{1}{2}} | \nabla f (z)|^2 [g_c (f ) (r, t)]^{- \frac{1}{2}} dm(t) \Big ) \right ]^{\frac{1}{2}} \\
&~~ \times \left [ \tau \Big ( \int^1_0 r \log \frac{1}{r} dr
\int_{\mathbb{T}} [g_c (f )(r, t)]^{\frac{1}{2}} | \nabla g (z)|^2 [g_c (f ) (r, t)]^{\frac{1}{2}} dm(t)  \Big ) \right ]^{\frac{1}{2}} \\
& = \frac{2}{\pi} \left [ \tau \Big ( \int^1_0 r \log \frac{1}{r} dr
\int_{\mathbb{T}} [g_c (f ) (r, t)]^{-1} | \nabla f (z)|^2 dm(t)  \Big ) \right ]^{\frac{1}{2}} \\
&~~ \times \left [ \tau \Big ( \int^1_0 r \log \frac{1}{r} dr
\int_{\mathbb{T}} g_c (f )(r, t) | \nabla g (z)|^2 d m(t)  \Big ) \right ]^{\frac{1}{2}} \\
& \triangleq \frac{2}{\pi} A \cdot B\end{split}\end{equation*}For
$A,$ since $r \log (1/r) \leq C(1-r^2)$ for $0<r<1$ we
have\begin{equation*}\begin{split}A^2 & \leq C \tau \Big (
\int_{\mathbb{T}} \int^1_0 [g_c (f ) (r, t)]^{-1} |\nabla f (r t)|^2  (1-r^2) d r d m (t) \Big )\\
& = C \tau \Big (  \int_{\mathbb{T}} \int^1_0 [g_c (f ) (t,r)]^{-1}
\frac{d g^2_c (f ) (t,r)}{dr} d r d m (t) \Big )\\
& = C \tau \Big ( \int_{\mathbb{T}} \int^1_0 \frac{d g_c (f )
(t,r)}{dr} d r d m (t) \Big )\\
& = C \| g_c (f ) \|_{L^1(L^{\infty}(\mathbb{T}) \otimes
\mathcal{M})}.\end{split}\end{equation*}

To estimate $B,$ we define\begin{equation*}D(j, k) = \{ (e^{2 \pi i
\theta}, r): (j-1) 2^{-k} \leq \theta < j 2^{-k}, ~~2^{-k-1 } < 1-r
\leq 2^{-k} \},\end{equation*}where $j =1,\ldots, 2^k$ and $k =0,
1,2,\cdots.$ Let $c_{jk} =(e^{2\pi i \theta_{jk}}, 1-2^{-k})$ with
$\theta_{jk} = (j-1/2)2^{-k}.$ Set\begin{equation*}[\tilde{A}_c(f,
\alpha)](e^{2\pi i \theta}, r) \triangleq [A_c(f, 3\pi
\alpha)](c_{jk}),~~\forall~(e^{2\pi i \theta}, r) \in D(j,
k),\end{equation*}and\begin{equation*}d_k (t) \triangleq
[\tilde{A}_c(f, \alpha)](t, 1 - 2^{-k-1}) - [\tilde{A}_c(f,
\alpha)](t, 1 - 2^{-k}) \geq 0,\end{equation*}respectively.
Since\begin{equation*}|z - e^{2\pi i \theta_{jk}}| \leq |z - e^{2\pi
i \theta}| + |e^{2\pi i \theta} - e^{2\pi i \theta_{jk}}| \leq
\alpha (1-|z|) + 2\pi (1 - |z|) \leq 3\pi \alpha
(1-|z|)\end{equation*}for $z \in \Gamma_{\alpha} (e^{2\pi i \theta},
1 - 2^{-k})$ with $(j-1) 2^{-k} \leq \theta < j 2^{-k},$ it is
concluded that\begin{equation*}[A_c(f, \alpha)](t, 1- 2^{-k}) \leq
[\tilde{A}_c(f, \alpha)](t, 1-2^{-k}) \leq [A_c(f, 5 \pi \alpha)](t,
1- 2^{-k})\end{equation*}for $t \in \mathbb{T}.$ Then, we have
\begin{equation*}\begin{split} B^2 & =\tau \Big ( \int^1_0
\int_{\mathbb{T}} g_c (f )
(t, r) |\nabla g (r t)|^2 r \log \frac{1}{r} d r d m(t) \Big )\\
& \leq C \tau \Big ( \int^1_0 \int_{\mathbb{T}} A_c (f ) (t, (1+r)/2
) |\nabla g (r t)|^2 r \log \frac{1}{r} d r d m(t) \Big )\\
& = C \tau \Big ( \int_{\mathbb{T}} d m (t)
\sum^{\infty}_{k=0}\int^{1 - 2^{-k-1}}_{1 -2^{-k}} A_c (f )
(t, (1+r)/2) | \nabla g (r t)|^2 r \log \frac{1}{r} d r \Big )\\
& \leq C \tau \Big ( \int_{\mathbb{T}} d m (t)
\sum^{\infty}_{k=0}[\tilde{A}_c(f, \alpha)](t, 1 - 2^{-k-2})
\int^{1 - 2^{-k-1}}_{1 -2^{-k}}  | \nabla g (r t)|^2 r \log \frac{1}{r} d r \Big )\\
& = C \tau \Big ( \int_{\mathbb{T}} d m (t) \sum^{\infty}_{k=0}\Big
( \sum^{k+1}_{j=0}d_j(t) \Big )
\int^{1 - 2^{-k-1}}_{1 -2^{-k}}  | \nabla g (r t)|^2 r \log \frac{1}{r} d r \Big )\\
& = C \tau \Big ( \int_{\mathbb{T}} d m (t) \Big ( d_0(t) \int^1_0 +
\sum^{\infty}_{k
= 1}d_k (t) \int^1_{1 -2^{-k+1}} \Big ) | \nabla g (r t)|^2 r \log \frac{1}{r} d r \Big )\\
& \leq C \tau[ d_0(1/2)] \Big \| \int_{\mathbb{T}}\int^1_0 | \nabla g (r t)|^2 r \log \frac{1}{r} d r d m (t) \Big \|\\
&~~+ C\sum^{\infty}_{k = 1} \sum^{2^k}_{j =1} \tau \big [ d_k ( e^{2
\pi \theta_{jk}}) \big ] \Big \| \int^{j2^{-k}}_{(j-1)2^{-k}} d
\theta \int^1_{1 -2^{-k+1}} | \nabla g (r e^{2\pi \theta})|^2 r \log
\frac{1}{r} d r \Big \|.\end{split}\end{equation*}Hence, by Theorem
\ref{th4.1} (2) one concludes that\begin{equation*}\begin{split} B^2
& \leq C \sum^{\infty}_{k = 0} \sum^{2^k}_{j =1}\tau \big [  d_k (
e^{2 \pi \theta_{jk}}) \big ] 2^{-k} \| g
\|^2_{\mathrm{BMO}_c}\\
& = C \| g \|^2_{\mathrm{BMO}_c}
\tau \Big ( \sum^{\infty}_{k = 0}\int_{\mathbb{T}} d_k ( t) d m(t) \Big )\\
& \leq C \| g \|^2_{\mathrm{BMO}_c} \tau \Big ( \int_{\mathbb{T}}
[A_c(f, 5 \pi \alpha)](t, 1) d m(t) \Big )\\
& \leq C \| g \|^2_{\mathrm{BMO}_c} \| A_c(f, 5 \pi \alpha)\|_{L^1}.
\end{split}\end{equation*}

Combining the estimates of $A$ and $B$ yields
that\begin{equation*}\begin{split}\|f\|_{\mathrm{H}^1_c} = \sup_{\|g
\|_{\mathrm{BMO}_c} \leq 1} \Big | \tau \Big ( \int_{\mathbb{T}} f
g^* d m \Big ) \Big | \leq C \| g_c (f, \alpha)
\|^{\frac{1}{2}}_{L^1} \| A_c(f, 5 \pi
\alpha)\|^{\frac{1}{2}}_{L^1}.\end{split}\end{equation*}Then, by
\eqref{6.2} we have that\begin{equation*}\begin{split} \| A_c(f, 5
\pi \alpha)\|_{L^1} \leq C \| g_c (f, \alpha)
\|_{L^1}.\end{split}\end{equation*}Therefore,\begin{equation*}\begin{split}\|f
\|_{\mathrm{H}^1_c} \leq C \| g_c (f, \alpha) \|_{L^1}.
\end{split}\end{equation*}This completes the proof.
\end{proof}

\begin{remark}\label{re6.1}
Let us consider the conformal map$$t (w) = i \frac{1-w}{1+w},\; w
\in \mathbb{T}, \; w \not= -1,
$$which maps $\mathbb{T}$ to the real line $\mathbb{R}.$ If $\varphi \in L^{\infty} (\mathcal{M},
L^2_c (\mathbb{R}, \frac{d t}{1+ t^2}))$ and if $\psi (w) = \varphi
(t (w)),$ then we see that $$\| \varphi \|_{\mathrm{BMO}_c
(\mathbb{R}, \mathcal{M})} \asymp \| \psi \|_{\mathrm{BMO}_c
(\mathbb{T}, \mathcal{M})},$$ following arguments in the classical
case \cite{G}. Consequently, by the $\mathrm{H}^1$-BMO duality we
have$$\| \varphi \|_{\mathrm{H}^1_c (\mathbb{R}, \mathcal{M})}
\asymp \| \psi \|_{\mathrm{H}^1_c (\mathbb{T}, \mathcal{M})}.$$Thus,
under the mapping $w \mapsto t (w),$ $\mathrm{H}^1_c (\mathbb{R},
\mathcal{M})$ and $\mathrm{H}^1_c (\mathbb{T}, \mathcal{M})$ are
transformed into each other.
\end{remark}

\begin{corollary}\label{cor6.1}
If $f \in \mathcal{H}^1 (\mathcal{M}, L^2_{cr} ( \mathbb{T}) ),$
then $f \in \mathcal{H}^1_{cr} (\mathbb{T}, \mathcal{M} )$ if and
only if there exist $g \in \mathcal{H}^1 (\mathcal{M}, L^2_c (
\mathbb{T}) )$ and $h \in \mathcal{H}^1 (\mathcal{M}, L^2_r (
\mathbb{T}) )$ such that $f = g+h$ with $A_c (g) \in L^1(\mathbb{T},
L^1(\mathcal{M}))$ and $A_r (h) \in L^1(\mathbb{T},
L^1(\mathcal{M})).$
\end{corollary}

\subsection*{Acknowledgments}
The author is grateful to Prof.Q.Xu for
useful suggestions and comments on the draft, and to Prof.T.N.Bekjan
for helpful discussions. I am also indebted to the anonymous referee
for useful comments and corrections.

\end{document}